\documentclass[10pt]{amsart}

\usepackage{amsthm,amsmath,amssymb}
\usepackage{textcomp,gensymb}
\usepackage[right=1.13in,left=1.13in,bottom=1in,top=1in]{geometry}
\usepackage{datetime}
\usepackage{multirow}
\usepackage{multicol}
\usepackage{graphicx}
\usepackage{booktabs}
\usepackage{caption}
\usepackage[numbers]{natbib}
\usepackage{listings}
\usepackage{makecell}
\usepackage{longtable}
\usepackage{float}
\usepackage{array}
\usepackage{titlesec}
\usepackage[T1]{fontenc} 

\newcolumntype{L}[1]{>{\raggedright\let\newline\\\arraybackslash\hspace{0pt}}m{#1}}
\newcolumntype{C}[1]{>{\centering\let\newline\\\arraybackslash\hspace{0pt}}m{#1}}
\newcolumntype{R}[1]{>{\raggedleft\let\newline\\\arraybackslash\hspace{0pt}}m{#1}}

\apptocmd{\thebibliography}{\raggedright}{}{}

\allowdisplaybreaks

\newcommand\Tstrut{\makebox[0pt][c]{\rule{0pt}{2.6ex}}}         
\def\mystrut(#1,#2){\vrule height #1pt depth #2pt width 0pt}

\makeatletter
\patchcmd{\@settitle}{\uppercasenonmath\@title}{}{}{}
\patchcmd{\@setauthors}{\MakeUppercase}{}{}{}
\patchcmd{\section}{\scshape}{}{}{}
\makeatother

\usepackage[table]{xcolor}
\definecolor{solv}{rgb}{0.65,1,0.65}
\definecolor{impo}{rgb}{1,0.65,1}
\definecolor{tent}{rgb}{1,1,0.65}
\definecolor{real}{rgb}{0.65,1,1}
\definecolor{unkn}{gray}{0.85}
\definecolor{newt}{rgb}{0.5,0.5,1}
\definecolor{na}{gray}{0.1}
\definecolor{bug}{rgb}{1,0.5,0.5}

\title{Some spherical coverings on S2 and their algebraic numbers}
\author{Randall L. Rathbun}
\email{randallrathbun@gmail.com}

\let\micro\micro
\let\perthousand\perthousand

\begin{document}

\begin{abstract}
Spherical coverings on the S2 sphere and their algebraic numbers are given for the putatively optimal global solutions for some n-congruent spherical caps with minimal radius to completely cover the S2 sphere. A few locally optimal solutions are also examined.
\end{abstract}

\maketitle

\setcounter{section}{0}
\raggedbottom

\section*{Optimal Spherical Coverings on the S2 Sphere}
\subsection*{Introduction}

The problem of covering the unit sphere S2 with $n$-congruent spherical caps such that they have minimal radius is a difficult problem. The figure \ref{fig:example} shown below, which has 12 congruent caps completely covering the S2 sphere, illustrates the situation. Please note how three caps typically intersect at a common point, this is a key to determining minimal solutions.
\begin{figure}[ht]
	\begin{center}
		\includegraphics[type=pdf,ext=pdf,read=pdf,height=2in,angle=0]{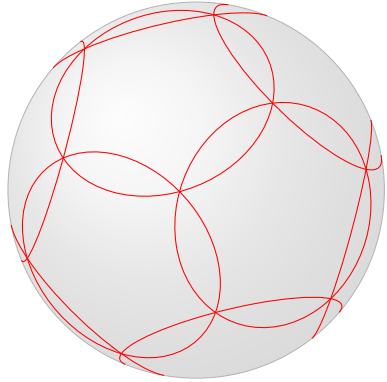}
		\caption{Spherical Caps example.}
		\label{fig:example}
	\end{center}
\end{figure}

\noindent\\[-2em]
It should be noted that a simple mathematical relationship between the cap radius and its distance from the centroid of the S2 sphere exists:
\begin{equation}
r = \sqrt{1-h^2}
\end{equation}
where $r$ is the radius of the caps and $h$ is the height of the cap center from the S2 sphere center. We seek to minimize $r$ or maximize $h$.

\noindent
The global solution is the minimal cap radius, maximal height, and their arrangement possible, for a set of $n$ spherical caps.

\section*{Putatively optimal global solutions}

Some putatively global solutions for $n\;=\;$ 2-10, 12, 14, 15, 16, 22, 32, \& 38 are given. Global solutions for $n=2-6$ and $n=10,\,12$ are already known.

The exact geometric structure has been successfully determined for four caps: $n=11, 13, 17, 20$; however no algebraic numbers were resolved, despite obtaining their spherical codes to hundreds or thousands of digits precision.

The two cases of 19 caps and 42 caps are discussed.

\subsection*{2 caps}
The optimal solution for 2 caps is trivial, since the whole sphere must be covered, the cap height is 0, thus leading to two hemispheres covering the sphere.
\begin{figure}[ht]
	\begin{center}
		\includegraphics[type=pdf,ext=pdf,read=pdf,height=2in,angle=0]{3.2.caps.}
		\caption{2 caps.}
		\label{fig:2caps}
	\end{center}
\end{figure}
\begin{center}
	\begin{tabular}{c|ccc}
		\multicolumn{4}{c}{Spherical code for 2 caps} \\
		\hline\Tstrut
		cap & x & y & z \\
		\hline\Tstrut
		1 & 1 & 0 & 0 \\
		2 & -1 & 0 & 0 \\
		\hline\Tstrut
		radius & 1 & & \\
		height & 0 & &
	\end{tabular}
\end{center}

\subsection*{3 caps}
The optimal solution for 3 caps is trivial, the 3 caps are arranged at the vertices of an equilateral triangle on an equator of the S2 sphere. Once again, the cap radius $r=1$ since the intersection at both poles must be a point to completely cover the sphere.
\begin{figure}[ht]
	\begin{center}
		\includegraphics[type=pdf,ext=pdf,read=pdf,height=2in,angle=0]{3.3.caps.}
		\caption{3 caps.}
		\label{fig:3caps}
	\end{center}
\end{figure}
\begin{center}
	\begin{tabular}{c|ccc}
		\multicolumn{4}{c}{Spherical code for 3 caps} \\
		\hline\Tstrut
		cap & x & y & z \\
		\hline\Tstrut
		1 & 1 & 0 & 0 \\
		2 & $-1/2$ & $\sqrt{3}/2$ & 0 \\
		3 & $-1/2$ & $-\sqrt{3}/2$ & 0 \\
		\hline\Tstrut
		radius & 1 & & \\
		height & 0 & &
	\end{tabular}
\end{center}

\subsection*{4 caps}
As might be surmised, the optimal solution for 4 caps is the tetrahedron, with the 4 caps at the vertices on the unit sphere.
\begin{figure}[ht]
	\begin{center}
		\includegraphics[type=pdf,ext=pdf,read=pdf,height=2in,angle=0]{3.4.caps.}
		\caption{4 caps.}
		\label{fig:4caps}
	\end{center}
\end{figure}
\begin{center}
	\begin{tabular}{c|ccc}
		\multicolumn{4}{c}{Spherical code for 4 caps} \\
		\hline\Tstrut
		cap & x & y & z \\
		\hline\Tstrut
		1 & 0 & 0 & 1 \\
		2 & $2\sqrt{2}/3$ & 0 & $-1/3$ \\
		3 & $-\sqrt{2}/3$ & $2/\sqrt{6}$ & $-1/3$ \\
		4 & $-\sqrt{2}/3$ & $-2/\sqrt{6}$ & $-1/3$ \\
		\hline\Tstrut
		radius & $2\sqrt{2}/3$ & & \\
		height & $1/3$ & &
	\end{tabular}
\end{center}

\subsection*{5 caps}
The optimal solution for 5 caps is the bi-pyramid, with the 3 caps at the vertices of an equilateral triangle on the equator and the two caps on the poles of the unit sphere.

\begin{figure}[H]
	\begin{center}
		\includegraphics[type=pdf,ext=pdf,read=pdf,height=2in,angle=0]{3.5.caps.}
		\caption{5 caps.}
		\label{fig:5caps}
	\end{center}
\end{figure}

\begin{center}
	\begin{tabular}{c|ccc}
		\multicolumn{4}{c}{Spherical code for 5 caps} \\
		\hline\Tstrut
		cap & x & y & z \\
		\hline\Tstrut
		1 & 0 & 0 & 1 \\
		2 & 1 & 0 & 0 \\
		3 & $-1/2$ & $\sqrt{3}/2$ & 0 \\
		4 & $-1/2$ & $-\sqrt{3}/2$ & 0 \\
		5 & 0 & 0 & -1 \\
		\hline\Tstrut \\[-1em]
		radius & $2/\sqrt{5}$ & & \\
		height & $1/\sqrt{5}$ & &
	\end{tabular}
\end{center}

\subsection*{6 caps}
The 6 caps for the optimal solution are on the vertices of the octahedron.

\begin{figure}[H]
	\begin{center}
		\includegraphics[type=pdf,ext=pdf,read=pdf,height=2in,angle=0]{3.6.caps.}
		\caption{6 caps.}
		\label{fig:6caps}
	\end{center}
\end{figure}

\begin{center}
	\begin{tabular}{c|ccc}
		\multicolumn{4}{c}{Spherical code for 6 caps} \\
		\hline\Tstrut
		cap & x & y & z \\
		\hline\Tstrut
		1 & 1 & 0 & 0 \\
		2 & -1 & 0 & 0 \\
		3 & 0 & 1 & 0 \\
		4 & 0 & -1 & 0 \\
		5 & 0 & 0 & 1 \\
		6 & 0 & 0 & -1 \\
		\hline\Tstrut\\[-1.1em]
		radius & $\sqrt{\frac{2}{3}}$ & & \\
		height & $1/\sqrt{3}$ & &
	\end{tabular}
\end{center}

\subsection*{7 caps}
The optimal solution for 7 caps are 5 caps on a regular pentagon on an equator and 2 caps at the poles.
\begin{figure}[H]
	\begin{center}
		\includegraphics[type=pdf,ext=pdf,read=pdf,height=2in,angle=0]{3.7.caps.}
		\caption{7 caps.}
		\label{fig:7caps}
	\end{center}
\end{figure}
\begin{center}
	\begin{tabular}{c|ccc}
		\multicolumn{4}{c}{Spherical code for 7 caps} \\
		\hline\Tstrut
		cap & x & y & z \\
		\hline\Tstrut
		1 & 0 & 0 & 1 \\
		2 & 0 & 0 & -1 \\
		3 & 0 & 1 & 0 \\
		4 & $\frac{-\sqrt{10+2\sqrt{5}}}{4}$ & $\frac{\sqrt{5}-1)}{4}$ & 0 \\[0.2em]
		5 & $\frac{-\sqrt(10-2\sqrt{5})}{4}$ & $\frac{-\sqrt{5}+1)}{4}$ & 0 \\[0.2em]
		6 & $\frac{\sqrt{10-2\sqrt{5}}}{4}$ & $\frac{-\sqrt{5}+1)}{4}$ & 0 \\[0.2em]
		7 & $\frac{\sqrt{10+2\sqrt{5}}}{4}$ & $\frac{\sqrt{5}-1)}{4}$ & 0 \\[0.2em]
		\hline\Tstrut\\[-1em]
		height & $\frac{\sqrt{7+2\sqrt{5}}}{\sqrt{29}}$ & & \\[0.5em]
		radius & $\frac{\sqrt{22-2\sqrt{5}}}{\sqrt{29}}$ & &
	\end{tabular}
\end{center}

\subsection*{8 caps}
The elegant solution for 8 caps places the caps at the vertices of 2 orthogonal tetrahedrons with their orientation reversed with respect to each other.
\begin{figure}[H]
	\begin{center}
		\includegraphics[type=pdf,ext=pdf,read=pdf,height=2in,angle=0]{3.8.caps.}
		\caption{8 caps.}
		\label{fig:8caps}
	\end{center}
\end{figure}
\noindent
Let
\begin{align*}
a & = 0.1715313637258024789\ldots = \text{a root of } 11x^{10} + 90x^8 + 192x^6 + 98x^4 - 139x^2 + 4 \\
b & = 0.9851786595630086290\ldots = \text{a root of } 11x^{10} - 145x^8 + 662x^6 - 1324x^4 + 1048x^2 - 256 \\
c & = 0.7867473620995775835\ldots = \text{a root of } 704x^{10} - 1281x^8 + 814x^6 - 160x^4 - 14x^2 + 1 \\
d & = 0.6172751317114242484\ldots = \text{a root of } 704x^{10} - 2239x^8 + 2730x^6 - 1636x^4 + 504x^2 - 64
\end{align*}
then the positions of the caps are:
\begin{longtable}[c]{c|rrr}
	\caption{Spherical code for 8 caps} \\
	cap & $x$ & $y$ & $z$ \\
	\hline\vspace*{-2.2ex}
	\endfirsthead
	\multicolumn{4}{c}%
	{\tablename\ \thetable\-- 8 caps code} \\
	cap & $x$ & $y$ & $z$ \\
	\hline\vspace*{-2.2ex}
	\endhead
	1 & $-a$ &  $b$ &  $0$ \\ 
	2 & $-a$ & $-b$ &  $0$ \\
	3 & $c$ &  $d$ &  $0$ \\
	4 & $a$ &  $0$ & $-b$ \\
	5 & $-c$ &  $0$ &  $d$ \\
	6 & $a$ &  $0$ &  $b$ \\
	7 & $-c$ &  $0$ & $-d$ \\
	8 & $c$ & $-d$ &  $0$
\end{longtable}
\noindent\\[-2em]
The height and radius are given by roots of the polynomial equations:
\begin{align*}
height & = 0.6673188865845566185\ldots = \text{a root of } 1331x^{10} - 2464x^8 + 1738x^6 - 1016x^4 + 291x^2 - 8 \\
radius & = 0.7447721152188417296\ldots = \text{a root of } 1331x^{10} - 4191x^8 + 5192x^6 - 2724x^4 + 272x^2 + 128 \\
\end{align*}

\subsection*{9 caps}
The solution for 9 caps places the caps at the vertices of 3 equilateral triangles, all co-planar, with one triangle at the equator, and the other 2 triangles equidistant above and below the equator. The smaller triangles are arranged exactly at $180\degree$ rotation with respect to the triangle on the equator.

\begin{figure}[H]
	\begin{center}
		\includegraphics[type=pdf,ext=pdf,read=pdf,height=2in,angle=0]{3.9.caps.}
		\caption{9 caps.}
		\label{fig:9caps}
	\end{center}
\end{figure}

\noindent
Let
\begin{align*}
a & = 0.6961773622954127151\ldots = \text{a root of } 16x^8 - 72x^6 + 129x^4 - 66x^2 + 9
\end{align*}
then the positions of the caps are:
\begin{longtable}[c]{c|ccr}
	\caption{Spherical code for 9 caps} \\
	cap & $x$ & $y$ & $z$ \\
	\hline\vspace*{-2.2ex}
	\endfirsthead
	\multicolumn{4}{c}%
	{\tablename\ \thetable\-- 9 caps code} \\
	cap & $x$ & $y$ & $z$ \\
	\hline\vspace*{-2.2ex}
	\endhead
	1 & $\sqrt{1-a^2}$ & 0 & $a$ \\[0.5em]
	2 & $\frac{-\sqrt{1-a^2}}{2}$ & $\frac{\sqrt{3}}{2}\sqrt{1-a^2}$ & $a$ \\[0.5em]
	3 & $\frac{-\sqrt{1-a^2}}{2}$ & $\frac{-\sqrt{3}}{2}\sqrt{1-a^2}$ & a \\[0.5em]
	4 & $1/2$ & $\frac{\sqrt{3}}{2}$ & 0 \\[0.5em]
	5 & $-1$ & 0 & 0 \\[0.5em]
	6 & $1/2$ & $\frac{-\sqrt{3}}{2}$ & 0 \\[0.5em]
	7 & $\sqrt{1-a^2}$ & 0 & $-a$ \\[0.5em]
	8 & $\frac{-\sqrt{1-a^2}}{2}$ & $\frac{\sqrt{3}}{2}\sqrt{1-a^2}$ & $-a$ \\[0.5em]
	9 & $\frac{-\sqrt{1-a^2}}{2}$ & $\frac{-\sqrt{3}}{2}\sqrt{1-a^2}$ & $-a$
\end{longtable}
\noindent\\[-2em]
The height and radius are given by roots of the two polynomial equations:
\begin{align*}
height & = 0.6961773622954127151\ldots = \text{a root of } 16x^8 - 72x^6 + 129x^4 - 66x^2 + 9 \\
radius & = 0.7178698212262454947\ldots = \text{a root of } 4x^4 - 4x^3 + 3x^2 + 4x - 4 \\
\end{align*}

\subsection*{10 caps}
The solution for 10 caps places 2 of the caps at the poles and the other 8 caps are arranged on 2 co-planar squares, rotated $90\degree$ with respect to each other, the squares equidistant on either side of the equator.

\begin{figure}[H]
	\begin{center}
		\includegraphics[type=pdf,ext=pdf,read=pdf,height=2in,angle=0]{3.10.caps.}
		\caption{10 caps.}
		\label{fig:10caps}
	\end{center}
\end{figure}
\noindent
The positions of the 10 caps are:
\begin{longtable}[c]{c|ccr}
	\caption{Spherical code for 10 caps} \\
	cap & $x$ & $y$ & $z$ \\
	\hline\vspace*{-2.2ex}
	\endfirsthead
	\multicolumn{4}{c}%
	{\tablename\ \thetable\-- 10 caps code} \\
	cap & $x$ & $y$ & $z$ \\
	\hline\vspace*{-2.2ex}
	\endhead
	1 & $0$ & $0$ & $1$ \\[0.5em]
	2 & $\sqrt{\sqrt{8}-2}$ & $0$ & $\sqrt{2}-1$ \\[0.5em]
	3 & $0$ & $\sqrt{\sqrt{8}-2}$ & $\sqrt{2}-1$ \\[0.5em]
	4 & $-\sqrt{\sqrt{8}-2}$ & $0$ & $\sqrt{2}-1$ \\[0.5em]
	5 & $0$ & $-\sqrt{\sqrt{8}-2}$ & $\sqrt{2}-1$ \\[0.5em]
	6 & $\sqrt{\sqrt{2}-1}$ & $\sqrt{\sqrt{2}-1}$ & $1-\sqrt{2}$ \\[0.5em]
	7 & $-\sqrt{\sqrt{2}-1}$ & $\sqrt{\sqrt{2}-1}$ & $1-\sqrt{2}$ \\[0.5em]
	8 & $-\sqrt{\sqrt{2}-1}$ & $-\sqrt{\sqrt{2}-1}$ & $1-\sqrt{2}$ \\[0.5em]
	9 & $\sqrt{\sqrt{2}-1}$ & $-\sqrt{\sqrt{2}-1}$ & $1-\sqrt{2}$ \\[0.5em]
	10 & $0$ & $0$ & $-1$
\end{longtable}
\noindent\\[-2em]
The height and radius are given by the two radicals:
\begin{equation*}
height = \sqrt{\frac{1+2\sqrt{2}}{7}} \qquad 
radius = \sqrt{\frac{6-2\sqrt{2}}{7}}
\end{equation*}

\subsection*{11 caps}
A first careful search for the algebraic numbers for 11 caps failed to resolve any using 96 digits accuracy. However the program was rerun to 655 digits accuracy and it turned out that a subtle geometric structure exists for this configuration discovered upon careful examination the second time.

The arrangement is a solitary cap at a pole and then 5 dipoles oriented along that polar axis in a balanced arrangement wrt the equatorial plane (see the spherical code below).

It takes 10 parameters $a-j$ to describe the caps layout.

\begin{figure}[H]
	\begin{center}
		\includegraphics[type=pdf,ext=pdf,read=pdf,height=2in,angle=0]{3.11.caps.}
		\caption{11 caps.}
		\label{fig:11caps}
	\end{center}
\end{figure}

\noindent
Let the approximate values (57 digits) for the 10 parameters, $a-j$, be set to:
\begin{align*}
a & = 0.634095032729788663209640739063305380363112029450554156101\ldots \\
b & = 0.770868938310986401236315482888284711235202039842114928182\ldots \\
c & = 0.336628849216948193713160388499986412982106827714811564941\ldots \\
d & = 0.149012786774920662896038564759947076002530372961658824364\ldots \\
e & = 0.121449107941445315463065360524596356873178461076451005551\ldots \\
f & = 0.796209636433271902876872055559046705042709312447749074551\ldots \\
g & = 0.389150175516029347544342728123496721953373440336695885810\ldots \\
a & = 0.768129234171283756456290915191894673622200146573358727739\ldots \\
i & = 0.842629988260650599505304915391079979121947505883657184842\ldots \\
j & = 0.159651353051667299154028178068364986420274751598506328331\ldots
\end{align*}

Then the spherical code for 11 caps is
\begin{longtable}[c]{c|ccr}
	\caption{Spherical code for 11 caps} \\
	cap & $x$ & $y$ & $z$ \\
	\hline\vspace*{-2.2ex}
	\endfirsthead
	\multicolumn{4}{c}%
	{\tablename\ \thetable\-- 11 caps code} \\
	cap & $x$ & $y$ & $z$ \\
	\hline\vspace*{-2.2ex}
	\endhead
	1 & $0$ & $0$ & $1$ \\
	2 & $b$ & $\sqrt{1-a^2-b^2}$ & $a$ \\
	3 & $-b$ & $-\sqrt{1-a^2-b^2}$ & $a$ \\
	4 & $-d$ & $\sqrt{1-c^2-d^2}$ & $c$ \\
	5 & $d$ & $-\sqrt{1-c^2-d^2}$ & $c$ \\
	6 & $f$ & $\sqrt{1-e^2-f^2}$ & $-e$ \\
	7 & $-f$ & $-\sqrt{1-e^2-f^2}$ & $-e$ \\
	8 & $h$ & $-\sqrt{1-g^2-h^2}$ & $-g$ \\
	9 & $-h$ & $\sqrt{1-g^2-h^2}$ & $-g$ \\
	10 & $j$ & $\sqrt{1-i^2-j^2}$ & $-i$ \\
	11 & $-j$ & $-\sqrt{1-i^2-j^2}$ & $-i$
\end{longtable}
\noindent\\[-2em]
The approximate height and radius (to 57 digits) are given by:
\begin{align*}
height & = 0.749797087499002078365989833146495343344798552462775132561\ldots \\
radius & = 0.661667837799309859442790840122286302673025957640146070376\ldots
\end{align*}
All 10 parameters and the height and radius have been determined to 2,022 digits, but no algebraic code was recovered for the global minimal arrangement of 11 caps. A local minimum has been recovered.

\subsection*{12 caps}
The solution for 12 caps places them at centers of the 12 faces of the regular dodecahedron which has 20 vertices and 30 edges. The solution has been previously determined as optimal, the dodecahedron known to the ancient Greeks.

\begin{figure}[H]
	\begin{center}
		\includegraphics[type=pdf,ext=pdf,read=pdf,height=2in,angle=0]{3.12.caps.}
		\caption{12 caps.}
		\label{fig:12caps}
	\end{center}
\end{figure}

\noindent
The Spherical code for 12 caps is
\begin{longtable}[c]{c|ccc}
	\caption{Spherical code for 12 caps} \\
	cap & $x$ & $y$ & $z$ \\
	\hline\vspace*{-2.2ex}
	\endfirsthead
	\multicolumn{4}{c}%
	{\tablename\ \thetable\-- 12 caps code} \\
	cap & $x$ & $y$ & $z$ \\
	\hline\vspace*{-2.2ex}
	\endhead
	1 & $1$ & $0$ & $0$ \\[0.5em]
	2 & $\sqrt{\frac{1}{5}}$ & $0$ & $2/\sqrt{5}$ \\[0.5em]
	3 & $\sqrt{\frac{1}{5}}$ & $\frac{-\sqrt{5+\sqrt{5}}}{10}$ & $\frac{5-\sqrt{5}}{10}$ \\[0.5em]
	4 & $\sqrt{\frac{1}{5}}$ & $\frac{-\sqrt{5-\sqrt{5}}}{10}$ & $\frac{-5-\sqrt{5}}{10}$ \\[0.5em]
	5 & $\sqrt{\frac{1}{5}}$ & $\frac{\sqrt{5-\sqrt{5}}}{10}$ & $\frac{-5-\sqrt{5}}{10}$ \\[0.5em]
	6 & $\sqrt{\frac{1}{5}}$ & $\frac{\sqrt{5+\sqrt{5}}}{10}$ & $\frac{5-\sqrt{5}}{10}$ \\[0.5em]
	7 & $-\sqrt{\frac{1}{5}}$ & $\frac{\sqrt{5-\sqrt{5}}}{10}$ & $\frac{5+\sqrt{5}}{10}$ \\[0.5em]
	8 & $-\sqrt{\frac{1}{5}}$ & $\frac{-\sqrt{5-\sqrt{5}}}{10}$ & $\frac{5+\sqrt{5}}{10}$ \\[0.5em]
	9 & $-\sqrt{\frac{1}{5}}$ & $\frac{-\sqrt{5+\sqrt{5}}}{10}$ & $\frac{-5+\sqrt{5}}{10}$ \\[0.5em]
	10 & $-\sqrt{\frac{1}{5}}$ & $0$ & $\frac{-2}{\sqrt{5}}$ \\[0.5em]
	11 & $-\sqrt{\frac{1}{5}}$ & $\frac{\sqrt{5+\sqrt{5}}}{10}$ & $\frac{-5+\sqrt{5}}{10}$ \\[0.5em]
	12 & $-1$ & $0$ & $0$
\end{longtable}
\noindent\\[-2em]
The height and radius are given by the two radicals
\begin{equation*}
height = \sqrt{\frac{5+2\sqrt{5}}{15}} \quad
radius = \sqrt{\frac{10-2\sqrt{5}}{15}}
\end{equation*}

\subsection*{13 caps}
The solution for 13 caps has successfully determined the geometric structures embedded into the cap placements, but even knowing the spherical code to 1,001 digits did not recover the algebraic numbers as hoped. The present spherical codes are a local minimum.

\begin{figure}[H]
	\begin{center}
		\includegraphics[type=pdf,ext=pdf,read=pdf,height=2in,angle=0]{3.13.caps.}
		\caption{13 caps.}
		\label{fig:13caps}
	\end{center}
\end{figure}
\vspace*{-2em}
\noindent
Let the approximate values (57 digits) for the 6 parameters, $a-f$, be set to:
\begin{align*}
a & = 0.866334167832381956179562484815574790082863963895771392911\ldots \\
b & = 0.845832418056934349532312578229913953813008606380422944999\ldots \\
c & = 0.524669634639377289135940108476275632025973178721753636316\ldots \\
d & = 0.0158150817488019292491381501149076972731869544937874573519\ldots \\
e & = 0.691828524265173209336956705953315739072959232420174267977\ldots \\
f & =  0.852938822579064225540289939438775325874868330870216691948\ldots
\end{align*}

Then the spherical code for 13 caps is
\begin{longtable}[c]{c|ccr}
	\caption{Spherical code for 13 caps} \\
	cap & $x$ & $y$ & $z$ \\
	\hline\vspace*{-2.2ex}
	\endfirsthead
	\multicolumn{4}{c}%
	{\tablename\ \thetable\-- 13 caps code} \\
	cap & $x$ & $y$ & $z$ \\
	\hline\vspace*{-2.2ex}
	\endhead
	1 & $\sqrt{1-a^2}$ & $0$ & $a$ \\
	2 & $-\sqrt{1-b^2}$ & $0$ & $b$ \\
	3 & $-d$ & $\sqrt{1-c^2-d^2}$ & $c$ \\
	4 & $-d$ & $-\sqrt{1-c^2-d^2}$ & $c$ \\
	5 & $1$ & $0$ & $0$ \\
	6 & $e$ & $\sqrt{1-e^2}$ & $0$ \\
	7 & $-f$ & $\sqrt{1-f^2}$ & $0$ \\
	8 & $-f$ & $-\sqrt{1-f^2}$ & $0$ \\
	9 & $e$ & $-\sqrt{1-e^2}$ & $0$ \\
	10 & $-d$ & $-\sqrt{1-c^2-d^2}$ & $-c$ \\
	11 & $-d$ & $\sqrt{1-c^2-d^2}$ & $-c$ \\
	12 & $-\sqrt{1-b^2}$ & $0$ & $-b$ \\
	13 & $\sqrt{1-a^2}$ & $0$ & $-a$
\end{longtable}
\noindent\\[-2em]
The approximate height and radius (to 57 digits) are given by:
\begin{align*}
height & = 0.797914989941050360712398793941776483164734262306857249527\ldots \\
radius & = 0.602769996621740871534109235592856793505094196214815579232\ldots
\end{align*}
All 6 parameters and height and radius have been determined to 1,001 digits, but no algebraic code was recovered for the global minimal arrangement of 13 caps. Most likely a local minimum has been found.

\subsection*{14 caps}
The algebraic solution for 14 caps has been successfully determined. There are 2 caps which appear at the poles, and the the 12 remaining caps are distributed on two co-planar hexagons, equally distributed on both sides of the equator.

\begin{figure}[H]
	\begin{center}
		\includegraphics[type=pdf,ext=pdf,read=pdf,height=2in,angle=0]{3.14.caps.}
		\caption{14 caps.}
		\label{fig:14caps}
	\end{center}
\end{figure}

\noindent
The the spherical code for 14 caps is
\begin{longtable}[c]{c|ccr}
	\caption{Spherical code for 14 caps} \\
	cap & $x$ & $y$ & $z$ \\
	\hline\vspace*{-2.2ex}
	\endfirsthead
	\multicolumn{4}{c}%
	{\tablename\ \thetable\-- 14 caps code} \\
	cap & $x$ & $y$ & $z$ \\
	\hline\vspace*{-2.2ex}
	\endhead
	1 & $0$ & $0$ & $1$ \\[0.7em]
	2 & $\sqrt{3\sqrt{3}-5}$ & $-\sqrt{9\sqrt{3}-15}$ & $2\,\sqrt{3}-3$ \\[0.7em]
	3 & $-\sqrt{3\sqrt{3}-5}$ & $\sqrt{9\sqrt{3}-15}$ & $2\,\sqrt{3}-3$ \\[0.7em]
	4 & $\sqrt{3\sqrt{3}-5}$ & $\sqrt{9\sqrt{3}-15}$ & $2\,\sqrt{3}-3$ \\[0.7em]
	5 & $-\sqrt{3\sqrt{3}-5}$ & $-\sqrt{9\sqrt{3}-15}$ & $2\,\sqrt{3}-3$ \\[0.7em]
	6 & $2\,\sqrt{3\sqrt{3}-5}$ & $0$ & $2\,\sqrt{3}-3$ \\[0.7em]
	7 & $-2\,\sqrt{3\sqrt{3}-5}$ & $0$ & $2\,\sqrt{3}-3$ \\[0.7em]
	8 & $0$ & $-2\,\sqrt{3\sqrt{3}-5}$ & $3-2\,\sqrt{3}$ \\[0.7em]
	9 & $0$ & $2\,\sqrt{3\sqrt{3}-5}$ & $3-2\,\sqrt{3}$ \\[0.7em]
	10 & $\sqrt{9\sqrt{3}-15}$ & $-\sqrt{3\sqrt{3}-5}$ & $3-2\,\sqrt{3}$ \\[0.7em]
	11 & $-\sqrt{9\sqrt{3}-15}$ & $\sqrt{3\sqrt{3}-5}$ & $3-2\,\sqrt{3}$ \\[0.7em]
	12 & $\sqrt{9\sqrt{3}-15}$ & $\sqrt{3\sqrt{3}-5}$ & $3-2\,\sqrt{3}$ \\[0.7em]
	13 & $-\sqrt{9\sqrt{3}-15}$ & $-\sqrt{3\sqrt{3}-5}$ & $3-2\,\sqrt{3}$ \\[0.7em]
	14 & $0$ & $0$ & $-1$
\end{longtable}

\noindent
The height and the radius are given by the radicals:
\begin{equation*}
height = \sqrt{\frac{6\sqrt{3}-3}{11}} \qquad
radius = \sqrt{\frac{6\sqrt{3}-14}{11}}
\end{equation*}

\subsection*{15 caps}
The algebraic solution for 15 caps has also been successfully determined. There are 6 parameters, used to construct the caps, and their algebraic numbers are known.

\begin{figure}[ht]
	\begin{center}
		\includegraphics[type=pdf,ext=pdf,read=pdf,height=2in,angle=0]{3.15.caps.}
		\caption{15 caps.}
		\label{fig:15caps}
	\end{center}
\end{figure}

\noindent
Let
\begin{align*}
a & = 0.7981658270508071459\ldots = \text{a root of } 866761x^{12} - 2571516x^{10} + 3920022x^8 - 3748572x^6 \\ & \qquad\qquad\qquad\qquad\qquad\qquad\qquad + 2360745x^4 - 956448x^2 + 186624 \\
b & = 0.4608212551057437875\ldots = \text{a root of } 931x^6 + 204x^5 - 438x^4 - 36x^3 + 123x^2 - 16 \\
c & = 0.4847751269560606253\ldots = \text{a root of } 14776336x^{12} - 27743256x^{10} + 23442345x^8 - 13991832x^6 \\ & \qquad\qquad\qquad\qquad\qquad\qquad\qquad + 6473520x^4 - 1772928x^2 + 186624 \\
d & = 0.2798850500445166073\ldots = \text{a root of } 3844x^6 + 2876x^5 - 127x^4 - 712x^3 - 196x^2 + 32x + 16 \\
e & = 0.9216425102114875750\ldots = \text{a root of } 931x^6 + 408x^5 - 1752x^4 - 288x^3 + 1968x^2 - 1024 \\
f & = 0.5597701000890332145\ldots = \text{a root of } 961x^6 + 1438x^5 - 127x^4 - 1424x^3 - 784x^2 + 256x + 256
\end{align*}
\noindent
Then the spherical code for 15 caps is
\begin{longtable}[c]{c|ccr}
	\caption{Spherical code for 15 caps} \\
	cap & $x$ & $y$ & $z$ \\
	\hline\vspace*{-2.2ex}
	\endfirsthead
	\multicolumn{4}{c}%
	{\tablename\ \thetable\-- 15 caps code} \\
	cap & $x$ & $y$ & $z$ \\
	\hline\vspace*{-2.2ex}
	\endhead
	1 & $-1/2$ & $0$ & $\sqrt{3}/2$ \\[0.3em]
	2 & $b$ & $\sqrt{1-a^2-b^2}$ & $a$ \\[0.3em]
	3 & $b$ & $-\sqrt{1-a^2-b^2}$ & $a$ \\[0.3em]
	4 & $-d$ & $\sqrt{1-c^2-d^2}$ & $c$ \\[0.3em]
	5 & $-d$ & $-\sqrt{1-c^2-d^2}$ & $c$ \\[0.3em]
	6 & $1$ & $0$ & $0$ \\[0.3em]
	7 & $-e$ & $\sqrt{1-e^2}$ & $0$ \\[0.3em]
	8 & $-e$ & $-\sqrt{1-e^2}$ & $0$ \\[0.3em]
	9 & $f$ & $\sqrt{1-f^2}$ & $0$ \\[0.3em]
	10 & $f$ & $-\sqrt{1-f^2}$ & $0$ \\[0.3em]
	11 & $-d$ & $-\sqrt{1-c^2-d^2}$ & $-c$ \\[0.3em]
	12 & $-d$ & $\sqrt{1-c^2-d^2}$ & $-c$ \\[0.3em]
	13 & $b$ & $-\sqrt{1-a^2-b^2}$ & $-a$ \\[0.3em]
	14 & $b$ & $\sqrt{1-a^2-b^2}$ & $-a$ \\[0.3em]
	15 & $-1/2$ & $0$ & $-\sqrt{3}/2$
\end{longtable}
\noindent
The height and the radius are given by a root of two polynomials:
\begin{align*}
height & = 0.8286479560382163503\ldots = \text{a root of } 923521x^{12} - 3229188x^{10} + 4897830x^8 - 3696996x^6 \\ & \qquad\qquad\qquad\qquad\qquad\qquad\qquad + 1421793x^4 - 272160x^2 + 20736 \\
radius & = 0.5597701000890332145\ldots = \text{a root of } 961x^6 + 1438x^5 - 127x^4 - 1424x^3 - 784x^2 \\ & \qquad\qquad\qquad\qquad\qquad\qquad\qquad + 256x + 256
\end{align*}

\subsection*{16 caps}
The algebraic solution for 16 caps has been determined.

\begin{figure}[ht]
	\begin{center}
		\includegraphics[type=pdf,ext=pdf,read=pdf,height=2in,angle=0]{3.16.caps.}
		\caption{16 caps.}
		\label{fig:16caps}
	\end{center}
\end{figure}

\noindent
The spherical code for 16 caps is
\begin{longtable}[c]{c|ccr}
	\caption{Spherical code for 16 caps} \\
	cap & $x$ & $y$ & $z$ \\
	\hline\vspace*{-2.2ex}
	\endfirsthead
	\multicolumn{4}{c}%
	{\tablename\ \thetable\-- 16 caps code} \\
	cap & $x$ & $y$ & $z$ \\
	\hline\vspace*{-2.2ex}
	\endhead
	1 & $0$ & $0$ & $1$ \\[0.7em]
	2 & $2\sqrt{\frac{2}{11}}$ & $0$ & $\sqrt{\frac{3}{11}}$ \\[0.7em]
	3 & $\sqrt{\frac{2}{11}}$ & $\sqrt{\frac{6}{11}}$ & $\sqrt{\frac{3}{11}}$ \\[0.7em]
	4 & $\sqrt{\frac{2}{11}}$ & $-\sqrt{\frac{6}{11}}$ & $\sqrt{\frac{3}{11}}$ \\[0.7em]
	5 & $-2\sqrt{\frac{2}{11}}$ & $0$ & $\sqrt{\frac{3}{11}}$ \\[0.7em]
	6 & $-\sqrt{\frac{2}{11}}$ & $-\sqrt{\frac{6}{11}}$ & $\sqrt{\frac{3}{11}}$ \\[0.7em]
	7 & $-\sqrt{\frac{2}{11}}$ & $\sqrt{\frac{6}{11}}$ & $\sqrt{\frac{3}{11}}$ \\[0.7em]
	8 & $-2\sqrt{\frac{2}{11}}$ & $-2\sqrt{\frac{2}{33}}$ & $-\sqrt{\frac{1}{33}}$ \\[0.7em]
	9 & $2\sqrt{\frac{2}{11}}$ & $-2\sqrt{\frac{2}{33}}$ & $-\sqrt{\frac{1}{33}}$ \\[0.7em]
	10 & $0$ & $4\sqrt{\frac{2}{33}}$ & $-\sqrt{\frac{1}{33}}$ \\[0.7em]
	11 & $0$ & $\frac{-2\sqrt{2}}{3}$ & $-\frac{1}{3}$ \\[0.7em]
	12 & $\sqrt{\frac{2}{3}}$ & $\frac{\sqrt{2}}{3}$ & $-\frac{1}{3}$ \\[0.7em]
	13 & $-\sqrt{\frac{2}{3}}$ & $\frac{\sqrt{2}}{3}$ & $-\frac{1}{3}$ \\[0.7em]
	14 & $0$ & $2\sqrt{\frac{2}{33}}$ & $-5\sqrt{\frac{1}{33}}$ \\[0.7em]
	15 & $\sqrt{\frac{2}{11}}$ & $-\sqrt{\frac{2}{33}}$ & $-5\sqrt{\frac{1}{33}}$ \\[0.7em]
	16 & $-\sqrt{\frac{2}{11}}$ & $-\sqrt{\frac{2}{33}}$ & $-5\sqrt{\frac{1}{33}}$
\end{longtable}
\noindent
The height and the radius are given by the two radicals:
\begin{equation*}
height=\sqrt{\frac{30+3\sqrt{33}}{67}} \qquad
radius=\sqrt{\frac{37-3\sqrt{33}}{67}}
\end{equation*}

\subsection*{17 caps}
Like 13 caps, the solution for 17 caps has successfully determined the spherical code to 501 digits by taking advantage of the embedded geometric structures to constrain the arrangement of caps, but unfortunately the algebraic numbers have not yet been recovered for this configuration. A local minima has been found instead.

\begin{figure}[ht]
	\begin{center}
		\includegraphics[type=pdf,ext=pdf,read=pdf,height=2in,angle=0]{3.17.caps.}
		\caption{17 caps.}
		\label{fig:17caps}
	\end{center}
\end{figure}

\noindent
Let the approximate values (57 digits) for the 8 parameters, $a-h$, be set to:
\begin{align*}
a & = 0.866081227467719654277982875756395580349279849736090829845\ldots \\
b & = 0.839492695431074414313697308819542902864844202373062181182\ldots \\
c & = 0.225436996700011596906179251271103687911086381469057624948\ldots \\
d & = 0.499142550440532090247633546193092719124405234631304652453\ldots \\
e & = 0.453462840795163529462601806038690222756460800405730272164\ldots \\
f & = 0.570638847472921471479793503369072549320101673699972073745\ldots \\
g & = 0.716685383316618478607701006566675126792988278713530663102\ldots \\
h & = 0.0723921505601956715934853706530802283110722765211476936016\ldots
\end{align*}

Then the spherical code for 17 caps using these 8 parameters is
\begin{longtable}[c]{c|ccr}
	\caption{Spherical code for 17 caps} \\
	cap & $x$ & $y$ & $z$ \\
	\hline\vspace*{-2.2ex}
	\endfirsthead
	\multicolumn{4}{c}%
	{\tablename\ \thetable\-- 17 caps code} \\
	cap & $x$ & $y$ & $z$ \\
	\hline\vspace*{-2.2ex}
	\endhead
	1 & $\sqrt{1-a^2}$ & $0$ & $a$ \\[0.5em]
	2 & $-c$ & $\sqrt{1-b^2-c^2}$ & $b$ \\[0.5em]
	3 & $-c$ & $-\sqrt{1-b^2-c^2}$ & $b$ \\[0.5em]
	4 & $-\sqrt{1-d^2}$ & $0$ & $d$ \\[0.5em]
	5 & $f$ & $\sqrt{1-e^2-f^2}$ & $e$ \\[0.5em]
	6 & $f$ & $-\sqrt{1-e^2-f^2}$ & $e$ \\[0.5em]
	7 & $1$ & $0$ & $0$ \\[0.5em]
	8 & $-g$ & $\sqrt{1-g^2}$ & $0$ \\[0.5em]
	9 & $-g$ & $-\sqrt{1-g^2}$ & $0$ \\[0.5em]
	10 & $-h$ & $\sqrt{1-h^2}$ & $0$ \\[0.5em]
	11 & $-h$ & $-\sqrt{1-h^2}$ & $0$ \\[0.5em]
	12 & $f$ & $-\sqrt{1-e^2-f^2}$ & $-e$ \\[0.5em]
	13 & $f$ & $\sqrt{1-e^2-f^2}$ & $-e$ \\[0.5em]
	14 & $-\sqrt{1-d^2}$ & $0$ & $-d$ \\[0.5em]
	15 & $-c$ & $-\sqrt{1-b^2-c^2}$ & $-b$ \\[0.5em]
	16 & $-c$ & $\sqrt{1-b^2-c^2}$ & $-b$ \\[0.5em]
	17 & $\sqrt{1-a^2}$ & $0$ & $-a$ 
\end{longtable}
\noindent\\[-2em]
The approximate height and radius (to 57 digits) are given by:
\begin{align*}
height & = 0.847187460907206953664488680584350129072560272215165104235\ldots \\
radius & = 0.531294086247531734784821102348202565559218797157927151241\ldots
\end{align*}
All 8 parameters and the height and radius have been determined to 501 digits, but no algebraic code has yet been recovered for these 17 caps on the S2 sphere.

\subsection*{18 caps}
After a preliminary run to 1,000 digits, a geometric structure was uncovered. The caps sit at the vertices of 6 co-planar sets of equilateral triangles, all aligned along a polar axis. 3 sets are above the equatorial plane, 3 below, and their distances from the equator are symmetrical.

There is a discrepancy in cap radius for this configuration, Tarnai \& G\'{a}sp\'{a}r give $30.013172\degree$ while the local minimum obtained here has a slighter larger value of $31.0132852\degree$, so this needs further investigation.

\begin{figure}[H]
	\begin{center}
		\includegraphics[type=pdf,ext=pdf,read=pdf,height=2in,angle=0]{3.18.caps.}
		\caption{18 caps.}
		\label{fig:18caps}
	\end{center}
\end{figure}
\noindent
Let the approximate values (38 digits) for the 8 parameters, $a-h$, be set to:
\begin{align*}
a & = 0.88869607772838093532058871775630035302\ldots \\
b & = 0.40536689571338645897257144124355299485\ldots \\
c & = 0.85189876064686144991449865810567729609\ldots \\
d & = 0.30143558361484884261468005967181211400\ldots \\
e & = 0.77266162426483545157943252333119354124\ldots \\
f & = 0.59487942737591607435684968399509107452\ldots \\
g & = 0.90957338356990010076719272668770552739\ldots \\
h & = 0.23290530593982136675359173767933586193\ldots
\end{align*}

The spherical code for 18 caps using these 8 parameters is
\begin{longtable}[c]{c|ccr}
	\caption{Spherical code for 18 caps} \\
	cap & $x$ & $y$ & $z$ \\
	\hline\vspace*{-2.2ex}
	\endfirsthead
	\multicolumn{4}{c}%
	{\tablename\ \thetable\-- 18 caps code} \\
	cap & $x$ & $y$ & $z$ \\
	\hline\vspace*{-2.2ex}
	\endhead
	1 & $\sqrt{\frac{1-a^2}{2}}$ & $\sqrt{\frac{1-a^2}{2}}$ & $a$ \\[0.7em]
	2 & $\frac{\left(\sqrt{6}-\sqrt{2}\right)\sqrt{1-a^2}}{4}$ & $-\frac{\left(\sqrt{6}+\sqrt{2}\right)\sqrt{1-a^2}}{4}$ & $a$ \\[0.7em]
	3 & $-\frac{\left(\sqrt{6}+\sqrt{2}\right)\sqrt{1-a^2}}{4}$ & $\frac{\left(\sqrt{6}-\sqrt{2}\right)\sqrt{1-a^2}}{4}$ & $a$ \\[0.7em]
	4 & $c$ & $-\sqrt{1-b^2-c^2}$ & $b$ \\[0.7em]
	5 & $-\frac{\sqrt{3}\sqrt{1-b^2-c^2}+c}{2}$ & $\frac{\sqrt{1-b^2-c^2}-c\sqrt{3}}{2}$ & $b$ \\[0.7em]
	6 & $\frac{\sqrt{3}\sqrt{1-b^2-c^2}-c}{2}$ & $\frac{\sqrt{1-b^2-c^2}+c\sqrt{3}}{2}$ & $b$ \\[0.7em]
	7 & $e$ & $\sqrt{1-d^2-e^2}$ & $d$ \\[0.7em]
	8 & $-\frac{\sqrt{3}\sqrt{1-d^2-e^2}+e}{2}$ & $-\frac{\sqrt{1-d^2-e^2}-e\sqrt{3}}{2}$ & $d$ \\[0.7em]
	9 & $\frac{\sqrt{3}\sqrt{1-d^2-e^2}-e}{2}$ & $-\frac{\sqrt{1-d^2-e^2}+e\sqrt{3}}{2}$ & $d$ \\[0.7em]
	10 & $f$ & $-\sqrt{1-d^2-f^2}$ & $-d$ \\[0.7em]
	11 & $-\frac{\sqrt{3}\sqrt{1-d^2-f^2}+f}{2}$ & $\frac{\sqrt{1-d^2-f^2}-f\sqrt{3}}{2}$ & $-d$ \\[0.7em]
	12 & $\frac{\sqrt{3}\sqrt{1-d^2-f^2}-f}{2}$ & $\frac{\sqrt{1-d^2-f^2}+f\sqrt{3}}{2}$ & $-d$ \\[0.7em]
	13 & $g$ & $\sqrt{1-b^2-g^2}$ & $-b$ \\[0.7em]
	14 & $-\frac{\sqrt{3}\sqrt{1-b^2-g^2}+g}{2}$ & $-\frac{\sqrt{1-b^2-g^2}-g\sqrt{3}}{2}$ & $-b$ \\[0.7em]
	15 & $\frac{\sqrt{3}\sqrt{1-b^2-g^2}-g}{2}$ & $-\frac{\sqrt{1-b^2-g^2}+g\sqrt{3}}{2}$ & $-b$ \\[0.7em]
	16 & $h$ & $\sqrt{1-a^2-h^2}$ & $-a$ \\[0.7em]
	17 & $-\frac{\sqrt{3}\sqrt{1-a^2-h^2}+h}{2}$ & $-\frac{\sqrt{1-a^2-h^2}-h\sqrt{3}}{2}$ & $-a$ \\[0.7em]
	18 & $\frac{\sqrt{3}\sqrt{1-a^2-h^2}-h}{2}$ & $-\frac{\sqrt{1-a^2-h^2}+h\sqrt{3}}{2}$ & $-a$
\end{longtable}
\noindent\\[-2em]
The approximate height and radius (to 38 digits) are given by:
\begin{align*}
height & = 0.85704785593771583133986734474450690675\ldots \\
radius & = 0.51523681218694408270171383057325677313\ldots
\end{align*}
All 8 parameters and the height and radius have been determined to 1,001 digits, but no algebraic code has yet been recovered for these 18 caps on the S2 sphere.

\subsection*{19 caps}
A careful search for the algebraic numbers for 19 caps failed to resolve any using 154 digits accuracy. Furthermore, none of the 19 caps lies upon any embedded polygons.

Like 11 caps, this set of caps which will require very high accuracy of the positions, in order to successfully recover the algebraic numbers for the spherical code.

Interestingly enough, the local minimum value found here for the cap radius, $30.374909\degree$ is less than that found by Tarnai \& G\'{a}sp\'{a}r of $30.382284\degree$. 
\begin{figure}[H]
	\begin{center}
		\includegraphics[type=pdf,ext=pdf,read=pdf,height=2in,angle=0]{3.19.caps.}
		\caption{19 caps.}
		\label{fig:19caps}
	\end{center}
\end{figure}

An approximate spherical code for an isometry of the 19 caps to 19 digits is:
\vspace*{2em}

\begin{longtable}[c]{c|rrr}
	\caption{Spherical code for 19 caps} \\
	cap & $x$ & $y$ & $z$ \\
	\hline\vspace*{-2.2ex}
	\endfirsthead
	\multicolumn{4}{c}%
	{\tablename\ \thetable\-- 19 caps code} \\
	cap & $x$ & $y$ & $z$ \\
	\hline\vspace*{-2.2ex}
	\endhead
	1 & 0 & 0 & 1 \\
	2 & 0 & 0.6349761541800811871 & 0.7725317363206988435 \\
	3 & 0.5054445214482143497 & -0.5118504076083359682 & 0.6946473896655524819 \\
	4 & -0.7581515192682984255 & 0.1596157978684256165 & 0.6322413074942168137 \\
	5 & -0.4234878603963508922 & -0.6656837752201766543 & 0.6144291200012692593 \\
	6 & 0.8047489411172299238 & 0.2607800469930158432 & 0.5332662645067824403 \\
	7 & 0.4723506300405221527 & 0.8667489431614193441 & 0.1600973198678991952 \\
	8 & -0.4912697990033500679 & 0.8582706651730744570 & 0.1484097365086817373 \\
	9 & 0.2729225847858185526 & -0.9535579251124018077 & 0.1274384014697249313 \\
	10 & -0.8942840794753883523 & -0.4466984223769217937 & -0.02676760435351778201 \\
	11 & 0.8984314716567872042 & -0.4346451142101566152 & -0.06248612189806215329 \\
	12 & -0.9586695281987397591 & 0.2018315013174024181 & -0.2005412196511457839 \\
	13 & -0.2926045265007563933 & -0.9045898152807171949 & -0.3099997696155048657 \\
	14 & 0.8863416803125210279 & 0.2644358420743626163 & -0.3800948712718964430 \\
	15 & 0.1494400093439340701 & 0.8646336895501477554 & -0.4796626590659143753 \\
	16 & 0.4329709803299691077 & -0.5421967394763222409 & -0.7201102873125412109 \\
	17 & -0.4664477027156708823 & 0.5036405554329888274 & -0.7271676089832547174 \\
	18 & -0.4539608950401538617 & -0.3346307578099949915 & -0.8257976517899835280 \\
	19 & 0.2994951172520864587 & 0.1824399201425274308 & -0.9364925788710486081
\end{longtable}
\noindent\\[-2em]
The height and radius found by direct search to 60 digits are:
\begin{align*}
height = 0.862735188616817442401043418772500763853821481432952011260923\ldots \\
radius = 0.505656003941715727509736496422975179334541566435250287496485\ldots
\end{align*}
NOTE: This height is less than the suggested height
\begin{align*}
suggested\_height = 0.8643750430693183107783654853864889893336\ldots = \frac{33441846558889}{38689046875000}
\end{align*}
given by the rational polynomial approximation function (see section) which means the more improvement could be possible for this local minimal solution.

\subsection*{20 caps}
The solution for 20 caps is similar to 17 caps and a geometric structure has been discovered, but the 1,040 digits found for the configuration are insufficient to recover the algebraic numbers. The arrangement consists of a hexagon at the equator, with caps on the vertices and 2 embedded rectangles for 8 caps. The other 6 caps are balanced with respect to the equator, 3 above and 3 below, equidistant along the normal axis to the hexagon.

\begin{figure}[!ht]
	\begin{center}
		\includegraphics[type=pdf,ext=pdf,read=pdf,height=2in,angle=0]{3.20.caps.}
		\caption{20 caps.}
		\label{fig:20caps}
	\end{center}
\end{figure}

\noindent
Let the approximate values (60 digits) for the 5 parameters, $a-e$, be set to:
\begin{align*}
a & = 0.984669201215740467462956375889364996551879460868097951710969\ldots \\
b & = 0.713862052440409939517513697308912610262081017787932024697946\ldots \\
c & = 0.661241128715695374063059099852736537986896799623001382774189\ldots \\
d & = 0.621807744273859137652915226525332113809858413969242767291768\ldots \\
e & = 0.427536691316927764829460356835178471781602631049905478364429\ldots
\end{align*}

The spherical code for 20 caps using these 5 parameters is
\begin{longtable}[c]{c|rcr}
	\caption{Spherical code for 20 caps} \\
	cap & $x$ & $y$ & $z$ \\
	\hline\vspace*{-2.2ex}
	\endfirsthead
	\multicolumn{4}{c}%
	{\tablename\ \thetable\-- 20 caps code} \\
	cap & $x$ & $y$ & $z$ \\
	\hline\vspace*{-2.2ex}
	\endhead
	1 & $0$ & $-\sqrt{1-a^2}$ & $a$ \\[0.5em]
	2 & $0$ & $\sqrt{1-b^2}$ & $b$ \\[0.5em]
	3 & $d$ & $\sqrt{1-c^2-d^2}$ & $c$ \\[0.5em]
	4 & $-d$ & $\sqrt{1-c^2-d^2}$ & $c$ \\[0.5em]
	5 & $c$ & $-\sqrt{1-c^2-d^2}$ & $d$ \\[0.5em]
	6 & $-c$ & $-\sqrt{1-c^2-d^2}$ & $d$ \\[0.5em]
	7 & $0$ & $-\sqrt{1-e^2}$ & $e$ \\[0.5em]
	8 & $a$ & $\sqrt{1-a^2}$ & $0$ \\[0.5em]
	9 & $-a$ & $\sqrt{1-a^2}$ & $0$ \\[0.5em]
	10 & $b$ & $-\sqrt{1-b^2}$ & $0$ \\[0.5em]
	11 & $-b$ & $-\sqrt{1-b^2}$ & $0$ \\[0.5em]
	12 & $e$ & $-\sqrt{1-e^2}$ & $0$ \\[0.5em]
	13 & $-e$ & $-\sqrt{1-e^2}$ & $0$ \\[0.5em]
	14 & $0$ & $-\sqrt{1-e^2}$ & $-e$ \\[0.5em]
	15 & $c$ & $-\sqrt{1-c^2-d^2}$ & $-d$ \\[0.5em]
	16 & $-c$ & $-\sqrt{1-c^2-d^2}$ & $-d$ \\[0.5em]
	17 & $d$ & $\sqrt{1-c^2-d^2}$ & $-c$ \\[0.5em]
	18 & $-d$ & $\sqrt{1-c^2-d^2}$ & $-c$ \\[0.5em]
	19 & $0$ & $\sqrt{1-b^2}$ & $-b$ \\[0.5em]
	20 & $0$ & $-\sqrt{1-a^2}$ & $-a$
\end{longtable}
\noindent\\[-2em]
The approximate height and radius (to 60 digits) are given by:
\begin{align*}
height & = 0.869295752152719094145304529712992668529115364736219859510985\ldots \\
radius & = 0.494292317651446385897084343851272189920942056249566920579577\ldots
\end{align*}
No algebraic code has yet been recovered for these 20 caps on the S2 sphere.

\newpage
\subsection*{Comment on next few caps listed}

The next few caps are those whose algebraic numbers were relatively easy to find. No attempt was made to follow the linear progression for $n$-caps for $n > 20$.

\subsection*{22 caps}
The algebraic numbers for 22 caps has been discovered. The arrangement of caps is successfully parameterized by 6 parameters, $a-f$, similar to 15 caps. There are 2 caps at the poles, and the other 20 caps are distributed on 4 co-planar irregular pentagons.

\begin{figure}[ht]
	\begin{center}
		\includegraphics[type=pdf,ext=pdf,read=pdf,height=2in,angle=0]{3.22.caps.}
		\caption{22 caps.}
		\label{fig:22caps}
	\end{center}
\end{figure}

\noindent
Let
\begin{align*}
a & = 0.6919090810112949212\ldots \text{a root of } 95x^6 + 30x^5 - 55x^4 - 20x^3 + 25x^2 - 10x - 1 \\
b & = 0.2231055241263849123\ldots \text{a root of } 9025x^{12} - 19550x^{10} + 15325x^8 - 5035x^6 + 615x^4 - 40x^2 + 1 \\
c & = 0.5840978452407353910\ldots \text{a root of } 9025x^{12} - 12550x^{10} + 8500x^8 - 2885x^6 + 465x^4 - 35x^2 + 1 \\
d & = 0.1677844402621232172\ldots \text{a root of } 25x^6 - 50x^5 + 75x^4 - 20x^3 - 25x^2 + 70x - 11 \\
e & = 0.3046362789712863538\ldots \text{a root of } 625x^{12} - 1250x^{10} + 1750x^8 - 925x^6 + 225x^4 - 25x^2 + 1 \\
f & = 0.7975481325531225265\ldots \text{a root of } 625x^{12} - 2500x^{10} + 3625x^8 - 1925x^6 + 375x^4 - 50x^2 + 1
\end{align*}
\noindent
Then the spherical code for 22 caps is
\begin{longtable}[c]{c|ccr}
	\caption{Spherical code for 22 caps} \\
	cap & $x$ & $y$ & $z$ \\
	\hline\vspace*{-2.2ex}
	\endfirsthead
	\multicolumn{4}{c}%
	{\tablename\ \thetable\-- Spherical code for 22 caps} \\
	cap & $x$ & $y$ & $z$ \\
	\hline\vspace*{-2.2ex}
	\endhead
	1 & $0$ & $0$ & $1$ \\[0.5em]
	2 & $\sqrt{1-a^2}$ & $0$ & $a$ \\[0.5em]
	3 & $b$ & $\sqrt{1-a^2-b^2}$ & $a$ \\[0.5em]
	4 & $-c$ & $\sqrt{1-a^2-c^2}$ & $a$ \\[0.5em]
	5 & $-c$ & $-\sqrt{1-a^2-c^2}$ & $a$ \\[0.5em]
	6 & $b$ & $-\sqrt{1-a^2-b^2}$ & $a$ \\[0.5em]
	7 & $-\sqrt{1-d^2}$ & $0$ & $d$ \\[0.5em]
	8 & $-e$ & $-\sqrt{1-d^2-e^2}$ & $d$ \\[0.5em]
	9 & $f$ & $-\sqrt{1-d^2-f^2}$ & $d$ \\[0.5em]
	10 & $f$ & $\sqrt{1-d^2-f^2}$ & $d$ \\[0.5em]
	11 & $-e$ & $\sqrt{1-d^2-e^2}$ & $d$ \\[0.5em]
	12 & $\sqrt{1-d^2}$ & $0$ & $-d$ \\[0.5em]
	13 & $e$ & $\sqrt{1-d^2-e^2}$ & $-d$ \\[0.5em]
	14 & $-f$ & $\sqrt{1-d^2-f^2}$ & $-d$ \\[0.5em]
	15 & $-f$ & $-\sqrt{1-d^2-f^2}$ & $-d$ \\[0.5em]
	16 & $e$ & $-\sqrt{1-d^2-e^2}$ & $-d$ \\[0.5em]
	17 & $-\sqrt{1-a^2}$ & $0$ & $-a$ \\[0.5em]
	18 & $-b$ & $-\sqrt{1-a^2-b^2}$ & $-a$ \\[0.5em]
	19 & $c$ & $-\sqrt{1-a^2-c^2}$ & $-a$ \\[0.5em]
	20 & $c$ & $\sqrt{1-a^2-c^2}$ & $-a$ \\[0.5em]
	21 & $-b$ & $\sqrt{1-a^2-b^2}$ & $-a$ \\[0.5em]
	22 & $0$ & $0$ & $-1$
\end{longtable}
\noindent
The height and the radius are given by a root of two polynomials:
\begin{align*}
height & = 0.8844990833734237006\ldots = \text{a root of } 28745x^{12} - 35070x^{10} + 12895x^8 - 2820x^6 \\ & \qquad\qquad\qquad\qquad\qquad\qquad\qquad + 375x^4 - 30x^2 + 1 \\
radius & = 0.4665419289962835370\ldots = \text{a root of } 28745x^{12} - 137400x^{10} + 268720x^8 - 272960x^6 \\ & \qquad\qquad\qquad\qquad\qquad\qquad\qquad + 149760x^4 - 40960x^2 + 4096
\end{align*}

\subsection*{32 caps}
It was discovered that 32 caps has embedded geometric polygons, indicating that the global solution is easily discoverable. This polyhedron has been noticed before as being highly symmetric.

\begin{figure}[H]
	\begin{center}
		\includegraphics[type=pdf,ext=pdf,read=pdf,height=2in,angle=0]{3.32.caps.}
		\caption{32 caps.}
		\label{fig:32caps}
	\end{center}
\end{figure}

\noindent
The spherical code for 32 caps is
\begin{longtable}[c]{c|ccr}
	\caption{Spherical code for 32 caps} \\
	cap & $x$ & $y$ & $z$ \\
	\hline\vspace*{-1.5ex}
	\endfirsthead
	\multicolumn{4}{c}%
	{\tablename\ \thetable\-- Spherical code for 32 caps} \\
	cap & $x$ & $y$ & $z$ \\
	\hline\vspace*{-1.5ex}
	\endhead
	1 & $\sqrt{\left(3-\sqrt{5}\right)/6}$ & $0$ & $\sqrt{\left(3+\sqrt{5}\right)/6}$ \\[0.5em] 
	2 & $-\sqrt{\left(3-\sqrt{5}\right)/6}$ & $0$ & $\sqrt{\left(3+\sqrt{5}\right)/6}$ \\[0.5em] 
	3 & $0$ & $\sqrt{\left(5-\sqrt{5}\right)/10}$ & $\sqrt{\left(5+\sqrt{5}\right)/10}$ \\[0.5em] 
	4 & $0$ & $-\sqrt{\left(5-\sqrt{5}\right)/10}$ & $\sqrt{\left(5+\sqrt{5}\right)/10}$ \\[0.5em] 
	5 & $\sqrt{1/3}$ & $\sqrt{1/3}$ & $\sqrt{1/3}$ \\[0.5em] 
	6 & $-\sqrt{1/3}$ & $\sqrt{1/3}$ & $\sqrt{1/3}$ \\[0.5em] 
	7 & $\sqrt{1/3}$ & $-\sqrt{1/3}$ & $\sqrt{1/3}$ \\[0.5em] 
	8 & $-\sqrt{1/3}$ & $-\sqrt{1/3}$ & $\sqrt{1/3}$ \\[0.5em] 
	9 & $\sqrt{\left(5+\sqrt{5}\right)/10}$ & $0$ & $\sqrt{\left(5-\sqrt{5}\right)/10}$ \\[0.5em] 
	10 & $-\sqrt{\left(5+\sqrt{5}\right)/10}$ & $0$ & $\sqrt{\left(5-\sqrt{5}\right)/10}$ \\[0.5em] 
	11 & $0$ & $\sqrt{\left(3+\sqrt{5}\right)/6}$ & $\sqrt{\left(3-\sqrt{5}\right)/6}$ \\[0.5em] 
	12 & $0$ & $-\sqrt{\left(3+\sqrt{5}\right)/6}$ & $\sqrt{\left(3-\sqrt{5}\right)/6}$ \\[0.5em] 
	13 & $-\sqrt{\left(3+\sqrt{5}\right)/6}$ & $\sqrt{\left(3-\sqrt{5}\right)/6}$ & $0$ \\[0.5em] 
	14 & $-\sqrt{\left(3+\sqrt{5}\right)/6}$ & $-\sqrt{\left(3-\sqrt{5}\right)/6}$ & $0$ \\[0.5em] 
	15 & $-\sqrt{\left(5-\sqrt{5}\right)/10}$ & $\sqrt{\left(5+\sqrt{5}\right)/10}$ & $0$ \\[0.5em] 
	16 & $-\sqrt{\left(5-\sqrt{5}\right)/10}$ & $-\sqrt{\left(5+\sqrt{5}\right)/10}$ & $0$ \\[0.5em] 
	17 & $\sqrt{\left(5-\sqrt{5}\right)/10}$ & $\sqrt{\left(5+\sqrt{5}\right)/10}$ & $0$ \\[0.5em] 
	18 & $\sqrt{\left(5-\sqrt{5}\right)/10}$ & $-\sqrt{\left(5+\sqrt{5}\right)/10}$ & $0$ \\[0.5em] 
	19 & $\sqrt{\left(3+\sqrt{5}\right)/6}$ & $\sqrt{\left(3-\sqrt{5}\right)/6}$ & $0$ \\[0.5em] 
	20 & $\sqrt{\left(3+\sqrt{5}\right)/6}$ & $-\sqrt{\left(3-\sqrt{5}\right)/6}$ & $0$ \\[0.5em] 
	21 & $0$ & $\sqrt{\left(3+\sqrt{5}\right)/6}$ & $-\sqrt{\left(3-\sqrt{5}\right)/6}$ \\[0.5em] 
	22 & $0$ & $-\sqrt{\left(3+\sqrt{5}\right)/6}$ & $-\sqrt{\left(3-\sqrt{5}\right)/6}$ \\[0.5em] 
	23 & $\sqrt{\left(5+\sqrt{5}\right)/10}$ & $0$ & $-\sqrt{\left(5-\sqrt{5}\right)/10}$ \\[0.5em] 
	24 & $-\sqrt{\left(5+\sqrt{5}\right)/10}$ & $0$ & $-\sqrt{\left(5-\sqrt{5}\right)/10}$ \\[0.5em] 
	25 & $\sqrt{1/3}$ & $\sqrt{1/3}$ & $-\sqrt{1/3}$ \\[0.5em] 
	26 & $-\sqrt{1/3}$ & $\sqrt{1/3}$ & $-\sqrt{1/3}$ \\[0.5em] 
	27 & $\sqrt{1/3}$ & $-\sqrt{1/3}$ & $-\sqrt{1/3}$ \\[0.5em] 
	28 & $-\sqrt{1/3}$ & $-\sqrt{1/3}$ & $-\sqrt{1/3}$ \\[0.5em] 
	29 & $0$ & $\sqrt{\left(5-\sqrt{5}\right)/10}$ & $-\sqrt{\left(5+\sqrt{5}\right)/10}$ \\[0.5em] 
	30 & $0$ & $-\sqrt{\left(5-\sqrt{5}\right)/10}$ & $-\sqrt{\left(5+\sqrt{5}\right)/10}$ \\[0.5em] 
	31 & $\sqrt{\left(3-\sqrt{5}\right)/6}$ & $0$ & $-\sqrt{\left(3+\sqrt{5}\right)/6}$ \\[0.5em] 
	32 & $-\sqrt{\left(3-\sqrt{5}\right)/6}$ & $0$ & $-\sqrt{\left(3+\sqrt{5}\right)/6}$ 
\end{longtable}
\noindent
The height and the radius are given by a root of two polynomials:
\begin{align*}
height & = 0.9226021945439894676\ldots = \text{a root of } 2245x^8 - 2480x^6 + 530x^4 - 40x^2 + 1 \\
radius & = 0.3857527584121915401\ldots = \text{a root of } 2245x^8 - 6500x^6 + 6560x^4 - 2560x^2 + 256
\end{align*}

\subsection*{38 caps}
The geometric arrangement of 38 caps also is highly regular, enabling the global solution to be found. The arrangement consists of two caps at the poles, and the other 36 caps are carefully distributed among 12 equilateral triangles, all co-planar, such that 6 are above and 6 below, and the distances of the triangle's centroids are symmetric with respect to their offsets from the sphere center.

\begin{figure}[ht]
	\begin{center}
		\includegraphics[type=pdf,ext=pdf,read=pdf,height=2in,angle=0]{3.38.caps.}
		\caption{38 caps.}
		\label{fig:38caps}
	\end{center}
\end{figure}

It was found that only 3 parameters, $a-c$, were needed to parameterize the triangles, and their algebraic numbers have been determined.

\noindent
Let
\begin{align*}
a & = 0.7996599850609895832\ldots = \text{a root of } 107x^{10} - 298x^9 + 1187x^8 - 312x^7 + 1110x^6 + 4356x^5 \\ & \qquad\qquad\qquad\qquad\qquad\qquad\qquad - 5346x^4 - 4536x^3 + 8991x^2 + 4374x - 6561 \\
b & = 0.4056080372332561015\ldots = \text{a root of } x^{10} - 62x^9 + 1453x^8 - 16408x^7 + 100418x^6 - 364132x^5 \\ & \qquad\qquad\qquad\qquad\qquad\qquad\qquad + 932866x^4 - 1625336x^3 + 1928317x^2 - 1408590x + 340881 \\
c & = 0.2282433142753381738\ldots = \text{a root of } 37x^{10} - 118x^9 + 685x^8 - 1460x^7 + 1490x^6 + 4816x^5 \\ & \qquad\qquad\qquad\qquad\qquad\qquad\qquad - 4766x^4 + 12596x^3 + 17305x^2 + 2022x - 1503
\end{align*}

\noindent
The spherical code for 38 caps is
\begin{longtable}[c]{c|ccr}
	\caption{Spherical code for 38 caps} \\
	cap & $x$ & $y$ & $z$ \\
	\hline\vspace*{-2.2ex}
	\endfirsthead
	\multicolumn{4}{c}%
	{\tablename\ \thetable\-- Spherical code for 38 caps} \\
	cap & $x$ & $y$ & $z$ \\
	\hline\vspace*{-2.2ex}
	\endhead
	1 & $0$ & $0$ & $1$ \\[0.5em]
	2 & $0$ & $\sqrt{1-a^2}$ & $a$ \\[0.5em]
	3 & $\sqrt{3}/2\sqrt{1-a^2}$ & $\sqrt{1-a^2}/2$ & $a$ \\[0.5em]
	4 & $\sqrt{3}/2\sqrt{1-a^2}$ & $-\sqrt{1-a^2}/2$ & $a$ \\[0.5em]
	5 & $0$ & $-\sqrt{1-a^2}$ & $a$ \\[0.5em]
	6 & $-\sqrt{3}/2\sqrt{1-a^2}$ & $-\sqrt{1-a^2}/2$ & $a$ \\[0.5em]
	7 & $-\sqrt{3}/2\sqrt{1-a^2}$ & $\sqrt{1-a^2}/2$ & $a$ \\[0.5em]
	8 & $\sqrt{1-b^2}$ & $0$ & $b$ \\[0.5em]
	9 & $\sqrt{1-b^2}/2$ & $\sqrt{3}/2\sqrt{1-b^2}$ & $b$ \\[0.5em]
	10 & $-\sqrt{1-b^2}/2$ & $\sqrt{3}/2\sqrt{1-b^2}$ & $b$ \\[0.5em]
	11 & $-\sqrt{1-b^2}$ & $0$ & $b$ \\[0.5em]
	12 & $-\sqrt{1-b^2}/2$ & $-\sqrt{3}/2\sqrt{1-b^2}$ & $b$ \\[0.5em]
	13 & $\sqrt{1-b^2}/2$ & $-\sqrt{3}/2\sqrt{1-b^2}$ & $b$ \\[0.5em]
	14 & $0$ & $\sqrt{1-c^2}$ & $c$ \\[0.5em]
	15 & $\sqrt{3}/2\sqrt{1-c^2}$ & $\sqrt{1-c^2}/2$ & $c$ \\[0.5em]
	16 & $\sqrt{3}/2\sqrt{1-c^2}$ & $-\sqrt{1-c^2}/2$ & $c$ \\[0.5em]
	17 & $0$ & $-\sqrt{1-c^2}$ & $c$ \\[0.5em]
	18 & $-\sqrt{3}/2\sqrt{1-c^2}$ & $-\sqrt{1-c^2}/2$ & $c$ \\[0.5em]
	19 & $-\sqrt{3}/2\sqrt{1-c^2}$ & $\sqrt{1-c^2}/2$ & $c$ \\[0.5em]
	20 & $\sqrt{1-c^2}$ & $0$ & $-c$ \\[0.5em]
	21 & $\sqrt{1-c^2}/2$ & $\sqrt{3}/2\sqrt{1-c^2}$ & $-c$ \\[0.5em]
	22 & $-\sqrt{1-c^2}/2$ & $\sqrt{3}/2\sqrt{1-c^2}$ & $-c$ \\[0.5em]
	23 & $-\sqrt{1-c^2}$ & $0$ & $-c$ \\[0.5em]
	24 & $-\sqrt{1-c^2}/2$ & $-\sqrt{3}/2\sqrt{1-c^2}$ & $-c$ \\[0.5em]
	25 & $\sqrt{1-c^2}/2$ & $-\sqrt{3}/2\sqrt{1-c^2}$ & $-c$ \\[0.5em]
	26 & $0$ & $\sqrt{1-b^2}$ & $-b$ \\[0.5em]
	27 & $\sqrt{3}/2\sqrt{1-b^2}$ & $\sqrt{1-b^2}/2$ & $-b$ \\[0.5em]
	28 & $\sqrt{3}/2\sqrt{1-b^2}$ & $-\sqrt{1-b^2}/2$ & $-b$ \\[0.5em]
	29 & $0$ & $-\sqrt{1-b^2}$ & $-b$ \\[0.5em]
	30 & $-\sqrt{3}/2\sqrt{1-b^2}$ & $-\sqrt{1-b^2}/2$ & $-b$ \\[0.5em]
	31 & $-\sqrt{3}/2\sqrt{1-b^2}$ & $\sqrt{1-b^2}/2$ & $-b$ \\[0.5em]
	32 & $\sqrt{1-a^2}$ & $0$ & $-a$ \\[0.5em]
	33 & $\sqrt{1-a^2}/2$ & $\sqrt{3}/2\sqrt{1-a^2}$ & $-a$ \\[0.5em]
	34 & $-\sqrt{1-a^2}/2$ & $\sqrt{3}/2\sqrt{1-a^2}$ & $-a$ \\[0.5em]
	35 & $-\sqrt{1-a^2}$ & $0$ & $-a$ \\[0.5em]
	36 & $-\sqrt{1-a^2}/2$ & $-\sqrt{3}/2\sqrt{1-a^2}$ & $-a$ \\[0.5em]
	37 & $\sqrt{1-a^2}/2$ & $-\sqrt{3}/2\sqrt{1-a^2}$ & $-a$ \\[0.5em]
	38 & $0$ & $0$ & $-1$
\end{longtable}
\noindent
The height and the radius are given by a root of two polynomials:
\begin{align*}
height & = 0.9331427893809348476\ldots \text{a root of } 2032489x^{20} - 8238102x^{18} + 22032105x^{16} - 37927368x^{14} \\ & \qquad\qquad\qquad\qquad\qquad\qquad\qquad + 43729362x^{12} - 36034308x^{10} + 20217114x^8 - 5697864x^6 \\ & \qquad\qquad\qquad\qquad\qquad\qquad\qquad + 312741x^4 - 144342x^2 - 19683 \\
radius & = 0.3595059590971591483\ldots \text{a root of } 2032489x^{20} - 12086788x^{18} + 39351192x^{16} - 85656480x^{14} \\ & \qquad\qquad\qquad\qquad\qquad\qquad\qquad + 129958848x^{12} - 137851392x^{10} + 99597312x^8 - 47652864x^6 \\ & \qquad\qquad\qquad\qquad\qquad\qquad\qquad + 15433728x^4 - 3407872x^2 + 262144
\end{align*}

\subsection*{42 caps}
The geometric arrangement of 42 caps appears to consist of 14 co-planar sets of equilateral triangles, with a triad of caps at the triangle vertices. These 14 triangles are arranged with 7 above and 7 below the equator, distributed symmetrically along the co-planar normal axis.

\begin{figure}[H]
	\begin{center}
		\includegraphics[type=pdf,ext=pdf,read=pdf,height=2in,angle=0]{3.42.caps.}
		\caption{42 caps.}
		\label{fig:42caps}
	\end{center}
\end{figure}

Several attempts using the SVD matrix factoring program to force the caps intersection points onto planes were unsuccessful in determining more than about 5 digits accuracy, and this was disappointing. The SVD matrix factoring program was used successfully for 13, 15 and 17 point sets, to discover the algebraic numbers (or spherical code to high accuracy). It was hoped that the same would occur here, but after running the algorithm for awhile, divergence occurred showing instability.

It is unclear how to proceed in further refining the caps positions in this elusive case.

\section*{Cross checking the minimal caps radius}
Table 2 of Tarnai and Gaspar's paper\cite{1} \textit{"Covering a sphere by equal circles, and the rigidity of its graph"} supplies the minimal cap radii that they found and was used to check the current results.

\begin{longtable}[c]{r|cc}
	\caption{Cross checking for minimal solutions} \\
	  & Tarnai - G\'{a}sp\'{a}r & Rathbun \\
	n & Radius$\degree$ & Radius$\degree$ \\
	\hline\vspace*{-1.5ex}
	\endfirsthead
	\multicolumn{3}{c}%
	{\tablename\ \thetable\-- Cross checking for minimal solutions} \\
	  & Tarnai - G\'{a}sp\'{a}r & Rathbun \\
	n & Radius$\degree$ & Radius$\degree$ \\
	\hline\vspace*{-1.5ex}
	\endhead
	2 & & 90.0000000000 \\
	3 & 90.000000 & 90.0000000000 \\
	4 & 70.528779 & 70.5287793655 \\
	5 & 63.434949 & 63.4349488229 \\
	6 & 54.735610 & 54.7356103172 \\
	7 & 51.026553 & 51.0265526631 \\
	8 & 48.138529 & 48.1395290861 \\
	9 & 45.878888 & 45.8788878287 \\
	10 & 42.307827 & 42.3078266301 \\
	11 & 41.427196 & 41.4271959586\makebox(0,6)[l]{\,\textdagger} \\
	12 & 37.377368 & 37.3773681406 \\
	13 & 37.068543 & 37.0685427025\makebox(0,6)[l]{\,\textdagger} \\
	14 & 34.937927 & 34.9379269231 \\
	15 & 34.039900 & 34.0399001237 \\
	16 & 32.898812 & 32.8988127601 \\
	17 & 32.092933 & 32.0929327861\makebox(0,6)[l]{\,\textdagger} \\
	18 & 31.013172 & 31.0132851551\makebox(0,6)[l]{\,\textasteriskcentered} \\
	19 & 30.382284 & 30.3749090533\makebox(0,6)[l]{\,$\diamond$} \\
	20 & 29.623096 & 29.6230957838\makebox(0,6)[l]{\,\textdagger} \\
	22 & 27.810059 & 27.8100587699 \\
	32 & 22.690480 & 22.6904803756 \\
	38 & 21.069858 & 21.0698583869 \\
	42 & 20.153842 & 20.2572026800\makebox(0,6)[l]{\,\textdagger} \\
\end{longtable}
\noindent
\textdagger \; -- definitely a local minimum, not globally optimal \\
\textasteriskcentered \; -- needs investigation, values differ slightly \\
$\diamond\;$ -- new lower value found

\section*{Algorithm used for processing the caps and common radius}

A simple GP-Pari algorithm was used for determining the caps planes and common height.

The caps normals are $V$ and the height is $h$ in the algorithm. The mesh command creates the spherical Voroni mesh for the $n$ caps. The height is determined by averaging the rotated caps from the mesh to the [0,0,1] orientation and taking the $z$ value which is actually the cap height $h$.

For a given cap, the intersection points with other caps all lie upon the same plane, thus created a polygon.

The SVD matrix factoring program finds the best solution from the first eigenvector and eigenvalues so all the cap intersection points for a given cap are as close as possible to the same plane (least squares fit). The normal to the plane is actually the normal of the caps vector. The final step is to parameterize the SVD results and re-obtain the caps positions $V$ again.

Running the algorithm several time converges the caps positions and height.

\noindent
\textbf{Caps and Height convergence algorithm - GP Pari code}

\noindent\texttt{
\begin{tabbing}
\{converge(V,h)= \\
\hphantom{00} \= local(prec,sz,m,avght,i,j,k,u,v,R,msz,p,W,c,K,U,pr,s,t); \\
\>	prec=default(realprecision); \\
\>	tny=1.0/10\^{}prec; \\
\>	sz=matsize(V)[2]; \\
\>	m=mesh(V,h); \\
\>	avght=0; \\
\>	for(i=1,sz, \\
\>	\hphantom{00} \= u=V[i];v=[0,0,1];R=rot3(u,v); \\
\>\>	msz=matsize(m[i])[2]; \\
\>\>	p=vector(msz,j,m[i][j]*R); \\
\>\>	avght+=sum(j=1,msz,p[j][3])/msz; \\
\>	); \\
\>	avht/=sz; \\
\>	W=V; \\
\>	for(k=1,sz, \\
\>\>	msz=matsize(m[k])[2]; \\
\>\>	c=centroid(m[k]); \\
\>\>	M=matrix(msz,3,i,j,precision(m[k][i][j]-c[j],prec)); \\
\>\>	K=matSVD(M,tny)[3]; \\
\>\>	U=unit\_vector(vector(3,i,K[i,1])); \\
\>\>	W[k]=U*sign(U[1])*sign(V[k][1]); \\
\>	); \\
\>	pr=P\_to\_parms(W); \\
\>	W=parms\_to\_P(pr); \\
\>	return([W,avght]); \\
\}
\end{tabbing}
}
\vspace*{-2em}
\section*{Future work}
The global minimum for 11 caps, 13 caps, 17 caps and 20 caps needs to be found, the minima given here are close, i.e. 1.0e-19, but only a local minima.

The difference for 18 caps in the cap radius needs to be carefully examined. During runs, it was noticed that convergence for the caps was chaotic, apparently the many local minima possibilities lock up the searching algorithm into a false global minimum.

The lack of convergence for 42 caps needs to be corrected. It is hoped to find the global solution for this symmetric arrangement.

\section*{Approximation function r(n) for caps radius size}

\begin{figure}[H]
	\begin{center}
		\includegraphics[width=5in]{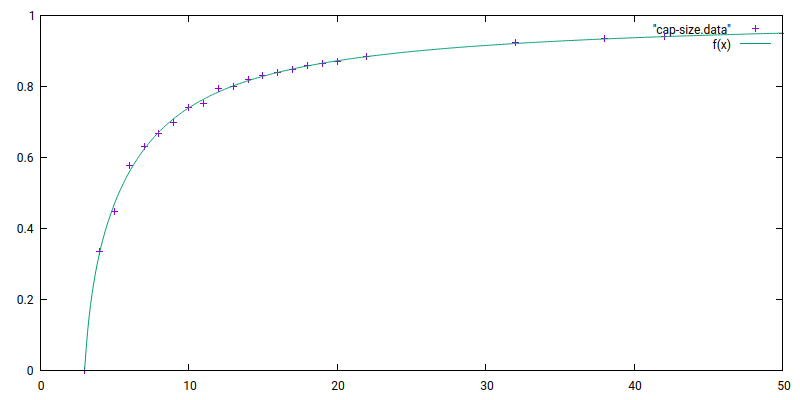}
		\caption{Caps radius size versus n.}
		\label{fig:capsrad}
	\end{center}
\end{figure}
The following rational polynomial $5^{th}$ degree approximation was found useful for determining the cap radius for a given $n$ for $ 3 <= n <= 150$ when initially starting searching and trial runs.

\begin{equation*}
r(n) = \frac{1.00185292n^5 - 2.67624769n^4 + 2.965272834n^3 - 43.610276n^2 + 217.5695441n - 366.876452}{n^5}
\end{equation*}

\section*{Acknowledgment}
Neil J.A. Sloane's caps coordinates were used as starting values for determining some of the caps solutions. See http://neilsloane.com/coverings/dim3/ for those spherical coverings.

\noindent
\textit{Randall L. Rathbun -- author \\ \today}\\[-3em]\hspace*{\fill}

\rule[1ex]{3.2in}{0.4pt}\hspace*{\fill} \\[-2.7em]\hspace*{\fill}

\end{document}